\documentclass[12pt,a4paper,russian,titlepage]{article}
\usepackage{babel,amsmath,amssymb,amsthm,graphicx,amsfonts,latexsym,epsfig,longtable,wrapfig}

\newtheorem{theorem}{Теорема}

\newtheorem{corollary}{Следствие}

\newcommand{\EHR}{{\mathrm{EHR}}}

\newcommand{\Pnp}{{\sf P}_{n, p}}

\newcommand{\D}{{\mathrm{D}}}

\newcommand{\K}{{\mathrm{K}}}
\newcommand{\MK}{{\mathrm{MK}}}
\newcommand{\NI}{{\mathrm{NI}}}

\begin{document}

\begin{center}{\Large Универсальный $k$--закон нуля или единицы}
\end{center}

\vspace{0.1cm}

\begin{center}{\large А.Д. Матушкин}
\end{center}

\vspace{0.5cm}

\section{Введение}
\label{intro}

В данной работе изучаются предельные вероятности свойств первого
порядка случайного графа в модели Эрдеша--Реньи
$G(n,n^{-\alpha})$, где $\alpha\in(0,1)	$. Мы нашли для любого $k$ и для любого рационального числа $t/s \in (0,1)$ интервал с правым концом $t/s$, на котором выполнен $k$--закон нуля или единицы, описывающий поведение вероятностей свойств первого порядка, выраженных формулами с ограниченной числом $k$ кванторной глубиной. Также для рациональных чисел $t/s$ с числителем, не превосходящим $2$, мы доказали, что логарифм длины найденного нами интервала имеет тот же порядок малости (при $n\rightarrow\infty$), что и логарифм длины наибольшего интеравла с правым концом $t/s$, на котором выполнен $k$--закон нуля или единицы.\\

Пусть $n\in\mathbb{N}$, $p=p(n)\in(0,1)$. Рассмотрим множество
$\Omega_n=\{G=(V_n,E)\}$ всех неориентированных графов без петель
и кратных ребер с множеством вершин $V_n=\{1,...,n\}$.
\emph{Случайный граф в модели Эрдеша--Реньи} (см.
\cite{Zhuk_Raigor}--\cite{Raigor}) это случайный элемент $G(n,p)$ со
значениями во множестве $\Omega_n$ и распределением ${\sf
P}_{n,p}$ на $(\Omega_n,\mathcal{F}_n)$, где
$\mathcal{F}_n=2^{\Omega_n},$ определенным формулой
$$
 {\sf P}_{n,p}(G)=p^{|E|}(1-p)^{C_n^2-|E|}, \quad\quad G\in\Omega_n.
$$

Говорят, что случайный граф {\it подчиняется закону нуля или
единицы} для класса свойств $\mathcal{C}$, если вероятность
выполнения каждого свойства из этого класса стремится либо к $0$,
либо к $1$.

Пожалуй, самым изученным в этом направлении является класс
свойств, выражаемых формулами первого порядка (см. \cite{Veresh},
\cite{Usp}). 
Класс свойств первого порядка обозначим $\mathcal{L}.$ В 1988 году Дж. Спенсер и С. Шела (см. \cite{Shelah},
\cite{Strange}) установили, что при $p=n^{-\alpha}$,
$\alpha\in(0, \infty)\setminus\mathbb{Q}$, случйный граф $G(n,p)$ подчиняется закону нуля или единицы для класса свойств первого порядка (в случае множества всех свойств первого порядка мы будем говорить просто ``закон нуля или единицы''). В той же работе было доказано, что при $\alpha\in(0,1]\cap\mathbb{Q}$ случайный граф $G(n,n^{-\alpha})$ закону нуля или единицы не подчиняется.

В данной работе мы рассмотрим класс $\mathcal{L}_k\subset \mathcal{L}$ свойств, выражаемых формулами с кванторной глубиной, ограниченной числом $k$ (см.~ \cite{Veresh}. Будем говорить, что случайный граф {\it подчиняется $k$--закону нуля или единицы}, если он подчиняется закону нуля или единицы для класса $\mathcal{L}_k$.\\ 

Ранее было установлено, что случайный граф $G(n,n^{-\alpha})$ подчиняется $k$--закону нуля или единицы при $\alpha$ близких к нулю и при $\alpha$ близких к единице.

\begin{theorem}[М.Е. Жуковский, \cite{Zhuk_law}, \cite{Zhuk_dan}]
Пусть $k\geq 3$ и $\alpha\in(0,1/(k-2))$. Тогда случайный граф
$G(n,n^{-\alpha})$ подчиняется $k$--закону нуля или единицы. Если
же $\alpha=1/(k-2)$, то случайный граф $G(n,n^{-\alpha})$ не
подчиняется $k$--закону нуля или единицы.
\label{closeToZero}
\end{theorem}

Тем самым, $\frac{1}{k-2}$ --- наменьшее положительное значение
$\alpha$, при котором нарушается $k$--закон нуля или единицы для
$G(n,n^{-\alpha}).$\\

\begin{theorem}[М.Е. Жуковский, \cite{Zhuk_extension_law}, \cite{Zhuk_large}]
Пусть $\beta$ --- произвольное положительное рациональное число,
$\alpha=1-\frac{1}{2^{k-1}+\beta}$. Пусть, кроме того,
$\mathcal{Q}$
--- множество положительных дробей с числителем, не превосходящим
числа $2^{k-1}$, $\widetilde{\mathcal{Q}}$ --- множество натуральных
чисел, не превосходящих $2^{k-1}-2$. Случайный граф
$G(n,n^{-{\alpha}})$ подчиняется $k$--закону нуля или единицы, если
$\beta\in(0,\infty)\setminus \mathcal{Q}$. Случайный граф
$G(n,n^{-{\alpha}})$ не подчиняется $k$--закону нуля или единицы,
если $\beta\in \widetilde{\mathcal{Q}}$. Если же
$\alpha\in\{1-\frac{1}{2^{k}-1},1-\frac{1}{2^{k}}\}$, то случайный
граф $G(n,n^{-{\alpha}})$ подчиняется $k$--закону нуля или единицы.
\label{closeToOne}
\end{theorem}

Итак, $1-\frac{1}{2^k-2}$ --- набольшее значение $\alpha$, меньшее
1, при котором нарушается $k$--закон нуля или единицы для
$G(n,n^{-\alpha}).$\\

Обратимся к точкам $\alpha$ из $(0,1)$, в окрестностях которых не выполнен $k$--закон нуля или единицы. Для этого введем понятие критической точки. Рациональная точка $\alpha$ называется {\it критической для свойства $A$}, если не выполнено следующее свойство. Существуют такие $\varepsilon>0$ и $\delta \in \{0,1\}$, что $\lim\limits_{n \rightarrow \infty}\Pnp(A) = \delta$ для всех $p \in [n^{-\alpha-\varepsilon}, n^{-\alpha+\varepsilon}]$. Дж. Спенсер (см. \cite{infiniteSpectra}) доказал, что существуют свойства первого порядка с бесконечным количеством критических точек в интервале $(0,1)$. В совместной работе Дж. Спенсер и С. Шела (см. \cite{criticalPoints}) доказали, что множество критических точек вполне упорядочено относительно порядка $>$. Таким образом, все его предельные  точки являются правосторонними пределами. Будем называть точку $k$--критической, если она является критической хотя бы для одного свойства первого порядка, кванторная глубина которого не превосходит числа $k$. Так как для любого натурального $k$ выполнено $|\mathcal{L}_k|<\infty$ (см.~\cite{Veresh}), то  слева от каждой рациональной точки $t/s\in(0,1)$ существует целый интервал, внутри которого нет ни одной $k$--критической точки (поэтому на этом интервале выполнен $k$--закон нуля или единицы). Мы нашли явно такой интервал для любых натуральных $k\geq4$ и рациональных $t/s\in(0,1)$, а также для рациональных чисел $t/s$ с числителем, не превосходящим $2$, мы доказали, что логарифм длины найденного нами интервала имеет тот же порядок малости, что и логарифм длины максимального такого интервала.

Дальнейшее повествование будет устроено следующим образом. В разделе~\ref{newResult} мы сформулируем доказанные нами теоремы о левых полуокрестностях рациональных чисел, внутри которых отсутствуют $k$--критические точки. Раздел~\ref{smallSubgraphs} посвящен предельным вероятностям свойств, связанных с содержанием некоторого фиксированного подграфа в случайном графе $G(n, n^{-\alpha})$. Данные свойства понадобятся нам при доказательстве теорем раздела~\ref{newResult}.

\section{Малые подграфы}
\label{smallSubgraphs}

Введем некоторые обозначения. Для произвольного графа $G$ обозначим за $v(G)$ число его вершин, $e(G)$ --- число его ребер, $a(G)$ --- число его автоморфизмов. Плотность графа $\frac{e(G)}{v(G)}$ обозначим $\rho(G)$. Далее мы эти обозначения будем использовать и в других разделах статьи. Граф $G$ называется \emph{сбалансированным}, если для каждого его подграфа $H$ выполнено неравенство $\rho(H)\leq \rho(G).$ Граф $G$ \emph{строго сбалансированный}, если для любого собственного подграфа $H\subset G$ справедливо строгое неравенство $\rho(H)<\rho(G).$ Сформулируем теорему (см. \cite{Erdos, Bol_small, Vince}) о количестве копий произвольного графа. Пусть $N_G$ --- количество копий $G$ в случайном графе $G(n,p)$. Положим
$\rho^{\max}(G)=\max\{\rho(H):\, H\subseteq G\}.$\\

\begin{theorem} [\cite{Erdos, Bol_small, Vince}] Пусть $G$ --- произвольный граф. Если $p=o\left(n^{-1/\rho^{\max}(G)}\right)$,
то
$$
 \lim_{n\rightarrow\infty}{\sf P}_{n,p}(N_G>0)=0.
$$
Если же $n^{-1/\rho^{\max}(G)}=o(p)$, то
$$
 \lim_{n\rightarrow\infty}{\sf P}_{n,p}(N_G>0)=1.
$$
Пусть теперь $G$
--- строго сбалансированный граф. Если $n^{-1/\rho(G)}=o(p),$ то
для любого $\varepsilon>0$ справедливо равенство
$$
 \lim\limits_{N\rightarrow\infty}{\sf P}_{n,p}(\left|N_G- {\sf
 E}_{n,p}N_G\right|\leq \varepsilon {\sf E}_{n,p}N_G)=1,
$$
где ${\sf E}_{n,p}$ --- это математическое ожидание по мере ${\sf
P}_{n,p}$. Если же $p=n^{-1/\rho(G)},$ то
$$
 \lim\limits_{n\rightarrow\infty}{\sf P}_{n,p}(N_G=0)=e^{-1/a(G)}.
$$
\label{small_subgraphs}
\label{erdos}
\end{theorem}

\begin{theorem} [A.Ruci\'{n}ski, A.Vince, \cite{graphExistence}] Пусть $\alpha = \frac{m}{n}$, $m,n\in\mathbb{N}$, $m<n$. Тогда существует строго сбалансированный граф с плотностью $\frac{1}{\alpha}$.
\label{strict_balanced_graph}
\end{theorem}

Так как свойство содержать фиксированный граф $H$ можно записать с помощью формулы первого порядка с конечной кванторной глубиной, то из теорем~\ref{erdos} и~\ref{strict_balanced_graph} следует, что для любого рационального $\alpha\in(0,1)$ существует $k$, для которого случайный граф $G(n,n^{-\alpha})$ не подчиняется $k$--закону нуля или единицы.

\section{Новые результаты}
\label{newResult}

Нашей целью будет исследовать выполнение закона нуля или единицы для $\alpha$, находящихся в малой левой полуокрестности некоторого фиксированного рационального числа. Ранее была доказана следующая теорема.

\begin{theorem} [М.Е. Жуковский, \cite{Zhuk_4}]
Если $\alpha\in(184/277,2/3)$, то случайный граф
$G(n,n^{-\alpha})$ подчиняется $4$--закону нуля или единицы.
\label{zhuk}
\end{theorem}

{\it Замечание.} Данный результат можно улучшить, заменив интервал $(184/277,2/3)$ на $(53/80,2/3)$, используя методы, аналогичные примененным в \cite{Zhuk_4}. Заметим также, что интервал $(184/277,2/3)$ является левой полуокрестностью точки $2/3$. Таким образом, данная теорема является аналогом следующих более общих фактов для $k=4$, $t/s = 2/3$.\\

\begin{theorem}
Пусть $\frac{t}{s}\in\mathbb{Q}\cap(0,1)$, $k\geq 4$. Положим $q=\frac{(s+1)^{k}-1}{s}$. Тогда внутри интервала $\left(\frac{tq}{sq+1},\frac{t}{s}\right)$ нет $k$--критических точек.
\label{main_interval}
\end{theorem}

\begin{corollary}
Пусть $\frac{t}{s}\in\mathbb{Q}\cap(0,1)$, $k\geq 4$. Положим $q=\frac{(s+1)^{k}-1}{s}$. Тогда если $\alpha\in\left(\frac{tq}{sq+1},\frac{t}{s}\right)$, то случайный граф $G(n,n^{-\alpha})$ подчиняется $k$--закону нуля или единицы.
\label{main_old}
\end{corollary}

Теорема~\ref{zhuk} дает больший интервал для $\alpha$, чем теорема~\ref{main_interval} для $k=4$, $t/s = 2/3$. Мы показали, что если выполнено $C_{k-2}^{[s/t]} < s+1$, то интервал, полученный в теореме~\ref{main_interval} можно увеличить, используя методы, примененные в доказательстве теоремы~\ref{zhuk}. Таким образом, теорему~\ref{main_interval} можно усилить, вследствие чего теорема~\ref{zhuk} будет частным случаем такого усиленного варианта.

\begin{theorem}
Пусть $\frac{t}{s}\in\mathbb{Q}\cap(0,1)$, $k\geq 4$. Положим $q=\frac{(s+1)^{k-2}(1+s(C^{[s/t]}_{k-2}+1))-1}{s}$. Тогда внутри интервала $\left(\frac{tq}{sq+1},\frac{t}{s}\right)$ нет $k$--критических точек.
\label{main_strong}
\end{theorem}

\begin{corollary}
Пусть $\frac{t}{s}\in\mathbb{Q}\cap(0,1)$, $k\geq 4$. Положим $q=\frac{(s+1)^{k-2}(1+s(C^{[s/t]}_{k-2}+1))-1}{s}$. Тогда если $\alpha\in\left(\frac{tq}{sq+1},\frac{t}{s}\right)$, то случайный граф $G(n,n^{-\alpha})$ подчиняется $k$--закону нуля или единицы.
\label{main_strong_old}
\end{corollary}

Легко видеть, что при $C_{k-2}^{[s/t]}<s+1$ теорема~\ref{main_strong} дает больший интервал значений $\alpha$, при которых выполнен $k$--закон, чем теорема~\ref{main_interval}.\\

Итак, ранее был доказан $k$--закон нуля или единицы для правой и левой полуокрестности 0 и 1 соответственно. При этом оставался огромный промежуток между $\frac{1}{k-2}$ и $1-\frac{1}{2^{k-1}}$, внутри которого не было найдено интервалов, для которых выполняется $k$--закон нуля или единицы. С помощью теоремы~\ref{main_interval} и теоремы~\ref{main_strong} мы можем найти такие интервалы на любом промежутке из интервала $(0,1)$.

Длина интервала, предъявленного в теореме~\ref{main_interval}, равна $\frac{t}{s(s+1)^{k}}$. Мы доказали, что для точек вида $\frac{t}{s}=\frac{2}{m}$, $m\geq 2$, теорема~\ref{main_interval} дает правильный порядок малости оценки логарифма длины максимального интервала с правым концом $\frac{t}{s}$, на котором нет критических точек. Более формально, для каждого $m\geq 2$, для каждого достаточно большого $k$ мы нашли слева от точки $\frac{2}{m}$ такую точку $\alpha$, для которой не выполнен $k$--закон нуля или единицы, причем величина $\frac{2}{m}-\alpha$ уменьшается экспоненциально с ростом $k$.

\begin{theorem}
Пусть $\frac{t}{s} = \frac{2}{m}$, $m\geq 2$, $k\geq 10m-5$. Тогда при $\alpha = \frac{t}{s} - \frac{1}{2^{k-10m+8}m(m-1)}$ случайный граф $G(n,n^{-\alpha})$ не подчиняется $k$--закону нуля или единицы.
\label{refutation}
\end{theorem}

Доказательство теорем~\ref{main_interval}, \ref{main_strong} и~\ref{refutation} будет построено по следующей схеме: в разделе~\ref{game} мы сформулируем теорему Эренфойхта, связывающую законы нуля или единицы и существование выигрышной стратегии
второго игрока в игре Эренфойхта. В разделе~\ref{proof_main_interval} мы
докажем теорему~\ref{main_interval}, в разделе~\ref{proof_main_strong} ---  теорему~\ref{main_strong}, а в разделе~\ref{proof_refutation} ---  теорему~\ref{refutation}. Раздел \ref{extensions} посвящен построению конструкций, свойства которых позволяют нам доказать, что игроки с необходимой вероятностью (стремящейся к 1) смогут играть в
соответствии с определенными нами в разделах~\ref{proof_main_interval} и~\ref{proof_main_strong} стратегиями.\\

\section{Расширения}
\label{extensions}

Обратимся к задаче, поставленной Дж.~Спенсером в 1990 году (см.
\cite{Alon}, \cite{Spencer}). Рассмотрим такие графы
$H,G,\widetilde{H},\widetilde{G}$, что $V(H)=\{x_1,...,x_k\}$,
$V(G)=\{x_1,...,x_l\}$,
$V(\widetilde{H})=\{\widetilde{x}_1,...,\widetilde{x}_k\}$,
$V(\widetilde{G})=\{\widetilde{x}_1,...,\widetilde{x}_l\}$, причем
$H\subset G$, $\widetilde{H}\subset\widetilde{G}$ (тем самым,
$k<l$). Граф $\widetilde{G}$ называется $(G,H)$--\emph{расширением
графа $\widetilde{H}$}, когда
$$
 \{x_{i_1},x_{i_2}\}\in E(G)\setminus E(H) \Rightarrow
 \{\widetilde{x}_{i_1},\widetilde{x}_{i_2}\}\in
 E(\widetilde{G})\setminus E(\widetilde{H}).
$$
Если выполняется соотношение
$$
 \{x_{i_1},x_{i_2}\}\in
 E(G)\setminus E(H) \Leftrightarrow
 \{\widetilde{x}_{i_1},\widetilde{x}_{i_2}\}\in
 E(\widetilde{G})\setminus E(\widetilde{H}),
$$
то $\widetilde{G}$ назовем \emph{точным $(G,H)$--расширением}.
Зафиксируем число $\alpha>0$. Положим
$$
 v(G,H)=|V(G)\setminus V(H)|, \,\,
 e(G,H)=|E(G)\setminus E(H)|,
$$
$$
 f_{\alpha}(G,H)=v(G,H)-\alpha e(G,H).
$$
Если для любого такого графа $S,$ что $H\subset S\subseteq G,$
выполнено неравенство $f_{\alpha}(S,H)>0,$ то пара $(G,H)$
называется \emph{$\alpha$--надежной} (см. \cite{Janson},
\cite{Alon}, \cite{Spencer}). Введем, наконец, понятие
максимальной пары. Пусть
$\widetilde{H}\subset\widetilde{G}\subset\Gamma$ и $T\subset K$,
причем $|V(T)|\leq|V(\widetilde{G})|.$ Пару
$(\widetilde{G},\widetilde{H})$ назовем
\emph{$(K,T)$--максимальной в $\Gamma$}, если у любого такого
подграфа $\widetilde{T}$ графа $\widetilde{G},$ что
$|V(\widetilde{T})|=|V(T)|$ и
$\widetilde{T}\cap\widetilde{H}\neq\widetilde{T},$ не существует
такого точного $(K,T)$--расширения $\widetilde{K}$ в
$\Gamma\setminus(\widetilde{G}\setminus\widetilde{T})$, что каждая
вершина из $V(\widetilde{K})\setminus V(\widetilde{T})$ не
соединена ребром ни с одной вершиной из $V(\widetilde{G})\setminus
V(\widetilde{T})$. Граф $\widetilde{G}$ назовем
\emph{$(K,T)$--максимальным в $\Gamma$}, если у любого такого
подграфа $\widetilde{T}$ графа $\widetilde{G},$ что
$|V(\widetilde{T})|=|V(T)|$, не существует такого точного
$(K,T)$--расширения $\widetilde{K}$ в
$\Gamma\setminus(\widetilde{G}\setminus\widetilde{T})$, что каждая
вершина из $V(\widetilde{K})\setminus V(\widetilde{T})$ не
соединена ребром ни с одной вершиной из $V(\widetilde{G})\setminus
V(\widetilde{T})$.

\begin{theorem} [J.H. Spencer, M.E. Zhukovskii, \cite{ZhukSpencer}]
Пусть $0 < \alpha_1 < \alpha_2 < 1$. Пусть пара $(G,H)$ --- $\alpha_2$--надежная (граф $H$ может быть пустым) и пусть $\mathcal{K}$ --- конечное множество пар графов, таких что для любой пары $(K,T) \in \mathcal{K}$ выполнены неравенства $f_{\alpha_1}(K,T) < 0$ и $v(T) \leq v(G)$. Пусть также $p\in [n^{-\alpha_2}, n^{-\alpha_1}]$. Тогда с вероятностью, стремящейся к единице при $n\rightarrow\infty$, выполнено следующее свойство. Для любых вершин $\widetilde{x}_1,...,\widetilde{x}_k$ графа $G(n,p)$ существует такое точное $(G,H)$--расширение $\widetilde{G}$ графа $\widetilde{H}:=G(n, p)|_{\{\widetilde{x}_1,...,\widetilde{x}_k\}}$ в $G(n,p)$, что пара $(\widetilde{G}, \widetilde{H})$ является $(K,T)$--максимальной в $G(n,p)$ для любой пары $(K,T) \in \mathcal{K}$.
\label{strong_maximal_graphs_rigid}
\end{theorem}

\section{Игра Эренфойхта}
\label{game}

Основным средством в доказательстве законов нуля или единицы для
свойств первого порядка случайных графов, как мы заметили выше,
служит теорема А. Эренфойхта, доказанная в 1960 году (см.
\cite{Ehren}). В данном разделе мы сформулируем ее частный случай
для графов. Прежде всего определим игру Эренфойхта $\EHR(G,H,i)$
на двух графах $G,H$ с количеством раундов, равным $i$ (см.,
например, \cite{Janson}, \cite{Alon}). Пусть
$V(G)=\{x_1,...,x_n\},$ $V(H)=\{y_1,...,y_m\}$. В $\nu\mbox{-}$ом
раунде ($1 \leq \nu \leq i$) Новатор выбирает вершину из любого
графа (он выбирает либо $x_{j_{\nu}}\in V(G)$, либо
$y_{j'_{\nu}}\in V(H)$). Затем Консерватор выбирает вершину из
оставшегося графа. Если Новатор выбирает в $\mu\mbox{-}$ом раунде,
скажем, вершину $x_{j_{\mu}}\in V(G),$ $j_{\mu}=j_{\nu}$
($\nu<\mu$), то Консерватор должен выбрать $y_{j'_{\nu}}\in V(H)$.
Если же в этом раунде Новатор выбирает, скажем, вершину
$x_{j_{\mu}}\in V(G),$ $j_{\mu}\notin\{j_1,...,j_{\mu-1}\},$ то и
Консерватор должен выбрать такую вершину $y_{j'_{\mu}}\in V(H)$,
что $j'_{\mu}\notin\{j'_1,...,j'_{\mu-1}\}.$ Если он не может
этого сделать, то игру выигрывает Новатор. К концу игры выбраны
вершины $x_{j_1},...,x_{j_{i}}\in V(G)$, а также вершины
$y_{j'_1},...,y_{j'_{i}}\in V(H)$. Некоторые из этих вершин могут
совпадать. Выберем из них только различные: $x_{h_1},...,x_{h_l};$
$y_{h'_1},...,y_{h'_l},$ $l \leq i.$ Консерватор побеждает тогда и
только тогда, когда соответствующие подграфы изоморфны с точностью
до порядка вершин:
$$
 G|_{\{x_{h_1},...,x_{h_l}\}}\cong
 H|_{\{y_{h'_1},...,y_{h'_l}\}}.
$$

\begin{theorem} [A. Ehrenfeucht, \cite{Ehren}]
Для любых двух графов $G,H$ и любого $i\in\mathbb{N}$ Консерватор
имеет выигрышную стратегию в игре $\EHR(G,H,i)$ тогда и только
тогда, когда для любого свойства $L$ первого порядка, выражаемого
формулой, кванторная глубина которой не превышает $i$, либо оба
графа обладают этим свойством, либо оба не обладают.
\end{theorem}

Несложно показать, что из этой теоремы вытекает следующее
следствие о законах нуля или единицы (см., например, \cite{ZhukSpencer}). Пусть ${\sf P}_{n,p(n)}\times {\sf
P}_{m,p(m)}$
--- прямое произведение мер ${\sf P}_{n,p(n)}, {\sf P}_{m,p(m)}$.

\begin{theorem} [J.H. Spencer, M.E. Zhukovskii, \cite{ZhukSpencer}]
Пусть $p_1 = p_1(n)$, $p_2 = p_2(n)$ --- две функции, все значения которых лежат в отрезке $[0,1]$, такие что $p_1(n) \leq p_2(n)$ для любого $n \in \mathbb{N}$. Тогда, если равенство
$$
 \lim\limits_{n,m\rightarrow\infty}{\sf P}_{n,p(n)}\times {\sf
 P}_{m,p(m)}(\{(A,B):\mbox{ у Консерватора есть выигрышная стратегия в игре}
$$
$$
 \EHR(A,B,k)\})=1
$$
выполнено для любой функции $p \in [p_1, p_2]$, то для любого свойства $L \in \mathcal{L}_k$ существует такое $\delta \in \{0,1\}$, что $\lim\limits_{n\rightarrow\infty}{\sf P}_{n, p(n)}\left(L\right) = \delta$ для любой функции $p \in [p_1, p_2]$.
\label{ehren_strong}
\end{theorem}

\section{Доказательство теоремы~\ref{main_interval}}
\label{proof_main_interval}

Пусть $\frac{t}{s}\in\mathbb{Q}\cap(0,1)$. Положим
\begin{equation}
q=\frac{(s+1)^{k}-1}{s}.
\label{q}
\end{equation}
Пусть $\alpha\in(\frac{tq}{sq+1},\frac{t}{s})$. Возьмем столь малое $\varepsilon>0$, чтобы отрезок $[\alpha-\varepsilon, \alpha+\varepsilon]$ полностью лежал внутри интервала $(\frac{tq}{sq+1},\frac{t}{s})$. Пусть $p\in[n^{-\alpha-\varepsilon}, n^{-\alpha+\varepsilon}]$.
Для каждого $i\in\{0,\dots,k\}$ положим
\begin{equation}
q_{i}=q-\frac{(s+1)^{k-i}-1}{s}.
\label{qi}
\end{equation}
Определим $\mathcal{K}_{i}$, $i\in\{0,\ldots,k\}$, как множество всех таких (попарно неизоморфных) пар $(K,T)$, что $f_{\alpha-\varepsilon}(K,T)<0$, $v(T)\leq q_{i}$, $v(K,T)\leq q-q_{i}$. Заметим, что $\mathcal{K}_0$ --- множество таких пар $(K,(\varnothing,\varnothing))$, что $\rho(K,T)>1/(\alpha-\varepsilon)$, а $\mathcal{K}_k=\varnothing$. Теперь определим множество графов $\mathcal{S}$. Граф $G$ принадлежит $\mathcal{S}$ тогда и только тогда, когда он обладает следующими свойствами.

\begin{itemize}
\item[{\sf 1)}]
В графе $G$ нет подграфов, количество вершин которых не превосходит $q$, а плотность больше, чем $\frac{1}{\alpha-\varepsilon}$.

\item[{\sf 2)}]

Пусть $\mathcal{H}$ --- множество таких $(\alpha+\varepsilon)$--надежных пар
$(H_{1},H_{2})$, что $v(H_{1})\leq q$. Тогда для любой пары
$(H_{1},H_{2})\in\mathcal{H}$ и для любого подграфа $G_{2}\subset G$ на $v(H_{2})$ вершинах в графе $G$ найдется такое точное
$(H_{1},H_{2})$--расширение $G_{1}$ подграфа $G_{2}$, что $(G_{1},G_{2})$ является $(K_{1},K_{2})$--максимальной парой в $G$ для любой пары $(K_{1},K_{2})\in\mathcal{K}_{i}$ для любого $i\in\{1,\ldots,k\}$.
\end{itemize}

В силу теорем~\ref{erdos} и~\ref{strong_maximal_graphs_rigid} для любого $p \in [n^{-\alpha-\varepsilon}, n^{-\alpha+\varepsilon}]$ справедливо $\lim\limits_{n \rightarrow \infty}\mathrm{P}_{n,p}\left(\mathcal{S}\right) = 1$. Следовательно, в силу теоремы~\ref{ehren_strong} для доказательства теоремы~\ref{main_interval} достаточно найти выигрышную стратегию Консерватора в игре $\mathrm{{EHR}}\left(G,H,k\right)$ для всех таких пар графов $\left(G,H\right)$, что $G,H\in\mathcal{S}$.

Итак, пусть $G,H\in\mathcal{S}$. Для каждого $h\in\{1,\dots,k\}$ обозначим за $X_{h}$, $Y_{h}$ графы выбранные в $h$-ом раунде Новатором и Консерватором соответственно. Положим также $X_0 = G$, $Y_0 = H$. Таким образом, множества $\{X_{h},Y_{h}\}$ и $\{G,H\}$ совпадают для всех $h\in\{0,\dots,k\}$. Вершины, выбранные в первых $h$ раундах в графе $X_{h}$, мы будем обозначать $x_{h}^{1},\dots,x_{h}^{h}$, а в графе $Y_{h}$ --- $y_{h}^{1},\dots,y_{h}^{h}$.\\

Определим теперь стратегию Консерватора для раундов $1,\dots,k$. Будем говорить, что {\it Консерватор выиграл в $l$-ом раунде}, где $l\in\{0,\ldots,k\}$, если существуют такие подграфы $\widetilde{X}_{l}\subset X_{l}$ и $\widetilde{Y}_{l}\subset Y_{l}$, что $\left\{ x_{l}^{1},\dots,x_{l}^{l}\right\} \in V(\widetilde{X}_{l})$, $\left\{ y_{l}^{1},\dots,y_{l}^{l}\right\} \in V(\widetilde{Y}_{l})$, и существует изоморфизм графов
$\varphi_{l}:V(\widetilde{X}_{l})\rightarrow V(\widetilde{Y}_{l})$, такой что $\varphi_{l}(x_{l}^{1})=y_{l}^{1},\dots,\varphi_{l}(x_{l}^{l})=y_{l}^{l}$, а также $v(\widetilde{X}_{l})=v(\widetilde{Y}_{l})\leq q_{l}$, и для любой пары $\left(K,T\right)\in\mathcal{K}_{l}$ графы $X_{l}$ и $Y_{l}$ не содержат $\left(K,T\right)$--расширений графов
$\widetilde{X}_{l}$ и $\widetilde{Y}_{l}$ соответственно. Заметим, что определение выигрыша Консерватора в игре $\EHR(G,H,k)$ совпадает с новым определением выигрыша в $k$-ом раунде. Заметим, кроме того, что для графов $\widetilde{X}_0 = (\varnothing, \varnothing)$, $\widetilde{Y}_0 = (\varnothing, \varnothing)$ в силу свойства~${\sf 1)}$ выполнено условие победы Консерватора в раунде с номером $0$. Докажем, что если $1\leq l\leq k$ и Консерватор выиграл в $(l-1)$-ом раунде, то он сможет выиграть и в $l$-ом раунде (в смысле нашего нового определения выигрыша), и, тем самым, докажем утверждение теоремы.

Итак, пусть было сыграно $l-1$ раундов, причем Консерватор выиграл в $(l-1)$-ом раунде. В $l$-ом раунде Новатор выбирает некоторую вершину $x_{l}^{l}$. Если $x_{l}^{l}\in V(\widetilde{X}_{l-1})$ или $x_{l}^{l}\in V(\widetilde{Y}_{l-1})$, то Консерватор, очевидно, выиграет в $l$-ом раунде, выбрав вершину $y_{l}^{l}$, являющуюся образом $x_{l}^{l}$ при изоморфизме $\varphi_{l-1}$ или $\varphi_{l-1}^{-1}$ соответственно.

Пусть $x_{l}^{l}\in V(X_{l-1})\setminus V(\widetilde{X}_{l-1})$ (в этом случае положим $\widetilde{X}_{l}=\widetilde{X}_{l-1}$, $\widetilde{Y}_{l}=\widetilde{Y}_{l-1}$) или $x_{l}^{l}\in V(Y_{l-1})\setminus V(\widetilde{Y}_{l-1})$ (в этом случае положим $\widetilde{X}_{l}=\widetilde{Y}_{l-1}$, $\widetilde{Y}_{l}=\widetilde{X}_{l-1})$.

Определим дополнительную конструкцию, которая нам потребуется в дальнейшем. Пусть $\mathcal{K}=\{(K,T),T\subset K\}$ --- некоторое конечное множество пар графов. Пару $(\widetilde K,\widetilde T)$ будем называть {\it $\mathcal{K}$--цепью}, если существует такое $r\in\mathbb{N}$
и такой набор графов $\widetilde T=K^0\subset K^1\subset\ldots\subset
K^r=\widetilde K$, что для любого $j\in\{1,\ldots,r\}$ граф $K^j$
является $(K,T)$--расширением графа $K^{j-1}$ для некоторой пары $(K,T)\in\mathcal{K}$.

Пусть пара $\left(Z,X_{l}|_{V(\widetilde{X}_{l})\cup\{x_{l}^{l}\}}\right)$ для некоторого $Z\subset X_{l}$ является такой $\mathcal{K}_{l}$--цепью, что граф $X_{l}$ не содержит $(K,T)$--расширений графа $Z$ ни при каком $(K,T)\in\mathcal{K}_{l}$. Докажем, что пара $(Z,\widetilde{X}_{l})$ является $(\alpha+\varepsilon)$--надежной. Пусть $X_{l}|_{V(\widetilde{X}_{l})\cup\{x_{l}^{l}\}}=K^{0}\subset K^{1}\subset\ldots\subset K^{r}=Z$, где для каждого $j\in\{1,\ldots,r\}$ граф $K^{j}$ является $(K,T)$--расширением графа $K^{j-1}$ для некоторой пары $(K,T)\in\mathcal{K}_{l}$. Для каждого $j\in\{1,\ldots,r\}$ положим $v_{j}=v(K^{j})-v(K^{j-1})$, $e_{j}=e(K^{j})-e(K^{j-1})$. В силу определения множества $\mathcal{K}_{l}$ для любого $j\in\{1,\ldots,r\}$ выполнено неравенство $f_{\alpha-\varepsilon}(K^{j},K^{j-1})<0$. Следовательно, $e_{j}>\frac{v_{j}}{\alpha-\varepsilon}>\frac{sv_{j}}{t}$. Таким образом, $e_{j}\geq\frac{sv_{j}+1}{t}$. Поэтому если $r\geq s+1$, то
\begin{eqnarray*}
f_{\alpha-\varepsilon}(K^r, \widetilde{X}_l) \leq
1 + \sum\limits_{j=1}^{r}v_j - (\alpha-\varepsilon)\sum\limits_{j=1}^{r}e_j \leq
1 - \dfrac{(\alpha-\varepsilon)r}{t} + \left(1-\dfrac{(\alpha-\varepsilon)s}{t} \right) \sum \limits_{j=1}^{r} v_j <
\\ < 1 - \dfrac{rq}{qs+1} + \dfrac{r(q-q_l)}{qs+1} =
1 - \dfrac{rq_l}{qs+1} \leq
1 - \dfrac{(s+1)q_1}{qs+1} =
\dfrac{s(q-q_1)+1-q_1}{qs+1} =
\\ = \dfrac{(s+1)^{k-1}-q_1}{qs+1} = 0
\end{eqnarray*}
в силу~(\ref{q})~и~(\ref{qi}). Итак, $f_{\alpha-\varepsilon}(K^{s+1},\widetilde{X}_{l})<0$. Так как $v(K^{s+1},\widetilde{X}_{l})\leq1+(s+1)(q-q_{l})=q-q_{l-1}$, получаем противоречие с определением графа $\widetilde{X}_{l}$. Следовательно, $r\leq s$. Заметим, что для любого подграфа $\widetilde{K}^{r} \subseteq K^{r}$, содержащего $\widetilde{X}_{l}$, выполнено неравенство $f_{\alpha-\varepsilon}(\widetilde{K}^{r},\widetilde{X}_{l}) \geq 0$, иначе аналогично получим противоречие с определением графа $\widetilde{X}_{l}$. Докажем, что выполнено более сильное неравенство $f_{\alpha+\varepsilon}(\widetilde{K}^{r},\widetilde{X}_{l}) > 0$. Действительно, иначе
$$
\frac{e(\widetilde{K}^{r},\widetilde{X}_{l})}{v(\widetilde{K}^{r},\widetilde{X}_{l})} \in \left[\frac{1}{\alpha+\varepsilon}, \frac{1}{\alpha-\varepsilon}\right] \subset \left(\frac{s}{t},\frac{sq+1}{tq}\right).
$$
При этом $v(\widetilde{K}^{r},\widetilde{X}_{l}) \leq q$, следовательно, $e(\widetilde{K}^{r},\widetilde{X}_{l}) \in \left(\frac{sq}{t},\frac{sq}{t} + \frac{1}{t}\right)$. Но данный интервал, очевидно, не содержит целых чисел. Таким образом, пара $(K^{r},\widetilde{X}_{l})$ действительно является $(\alpha+\varepsilon)$--надежной. Так как $v(K^{r})\leq q_{l-1}+1+\sum_{j=1}^{r}v_{j}\leq q_{l-1}+1+s\cdot(q-q_{l})=q_{l}$, то по свойству 2) в графе $Y_{l}$ существует такая вершина $y_{l}^{l}$ (которую и выберет Консерватор) и подграф $W$, являющийся точным $(K^{r},\widetilde{X}_{l})$--расширением графа $\widetilde{Y}_{l}$, что вершина $y_{l}^{l}$ является образом вершины $x_{l}^{l}$ при изоморфизме графов $\varphi_{l}:V(K^{r})\rightarrow V(W)$, таком что $\varphi_{l}|_{V(\widetilde{X}_{l})}=\varphi_{l-1}$ (или $\varphi_{l-1}^{-1}$). Переобозначим $\widetilde{X}_{l} = K^{r}$,
$\widetilde{Y}_{l} =  W$. Итак, Консерватор выиграл в $l$-ом раунде.

Таким образом, действуя по описанной стратегии, Консерватор выиграет в $k$-ом раунде. А значит, он выиграет в игре $\mathrm{EHR}\left(G,H,k\right)$, что и требовалось доказать.

\section{Доказательство теоремы~\ref{main_strong}}
\label{proof_main_strong}

Пусть $\frac{t}{s}\in\mathbb{Q}\cap(0,1)$. Положим $m=C^{[s/t]}_{k-2}+1$, $q=\frac{(s+1)^{k-2}(1+sm)-1}{s}$. Пусть $\alpha\in(\frac{tq}{sq+1},\frac{t}{s})$, $p~=~n^{-\alpha}$.
Введем некоторые дополнительные обозначения.
Для каждого $i\in\{0,\dots, k-2\}$ положим $q_{i}=q-\frac{(s+1)^{k-2-i}(1+sm)-1}{s}$. Как и при доказательстве теоремы~\ref{main_interval} определим $\mathcal{K}_{i}$ (но теперь для $i\in\{0,\ldots,k-2\}$) как множество всех таких пар $(K,T)$, что $f_{\alpha-\varepsilon}(K,T)<0$, $v(T)\leq q_{i}$, $v(K,T)\leq q-q_{i}$. Множество графов, обладающих свойствами ${\sf 1)}$ и ${\sf 2)}$, определяется в точности, как и при доказательстве теоремы~\ref{main_interval} (с новыми $\mathcal{K}_i$, где $i\in\{0,\ldots,k-2\}$, и $q$).

В силу теорем~\ref{erdos}, ~\ref{strong_maximal_graphs_rigid} и~\ref{ehren_strong} для доказательства теоремы~\ref{main_strong} достаточно найти выигрышную стратегию Консерватора в игре $\mathrm{{EHR}}\left(G,H,k\right)$ для всех таких пар графов $\left(G,H\right)$, что $G,H\in\mathcal{S}$.

Итак, пусть $G,H\in\mathcal{S}$. Как и ранее, для каждого $h\in\{1,\dots,k\}$ обозначим за $X_{h}$, $Y_{h}$ графы выбранные в $h$-ом раунде Новатором и Консерватором соответственно. Вершины, выбранные в первых $h$ раундах в графе $X_{h}$, мы будем обозначать $x_{h}^{1},\dots,x_{h}^{h}$, а в графе $Y_{h}$ --- $y_{h}^{1},\dots,y_{h}^{h}$.\\

Стратегия Консерватора для раундов $1,\dots,k-2$ в точности повторяет стратегию Консерватора для раундов $1,\dots,k$ в доказательстве теоремы~\ref{main_interval}. Как и ранее, мы говорим, что {\it Консерватор выиграл в $l$-ом раунде} (но теперь $l\in\{0,\ldots,k-2\}$), если существуют такие подграфы $\widetilde{X}_{l}\subset X_{l}$ и $\widetilde{Y}_{l}\subset Y_{l}$, что $\left\{ x_{l}^{1},\dots,x_{l}^{l}\right\} \in V(\widetilde{X}_{l})$, $\left\{ y_{l}^{1},\dots,y_{l}^{l}\right\} \in V(\widetilde{Y}_{l})$, и существует изоморфизм графов
$\varphi_{l}:V(\widetilde{X}_{l})\rightarrow V(\widetilde{Y}_{l})$, такой что $\varphi_{l}(x_{l}^{1})=y_{l}^{1},\dots,\varphi_{l}(x_{l}^{l})=y_{l}^{l}$, а также $v(\widetilde{X}_{l})=v(\widetilde{Y}_{l})\leq q_{l}$ и для любой пары $\left(K,T\right)\in\mathcal{K}_{l}$ графы $X_{l}$ и $Y_{l}$ не содержат $\left(K,T\right)$--расширений графов
$\widetilde{X}_{l}$ и $\widetilde{Y}_{l}$ соответственно. Действуя по описанной в разделе~\ref{proof_main_interval} стратегии, Консерватор выиграет в раунде с номером $k-2$.

Если в $(k-1)$-ом раунде Новатор выбирает вершину $x_{k-1}^{k-1}$ в графе $\widetilde X_{k-2}$ (в этом случае положим $\widetilde X_{k-1}=\widetilde X_{k-2}$, $\widetilde Y_{k-1}=\widetilde Y_{k-2}$) или в графе $\widetilde Y_{k-2}$ (в этом случае положим $\widetilde X_{k-1}=\widetilde Y_{k-2}$, $\widetilde Y_{k-1}=\widetilde X_{k-2}$), то Консерватор выберет вершину $y_{k-1}^{k-1}$, являющуюся образом вершины $x_{k-1}^{k-1}$ при изоморфизме $\varphi_{k-2}$ или $\varphi_{k-2}^{-1}$ соответственно. Если и в $k$-ом раунде Новатор выберет вершину либо в графе $\widetilde X_{k-1}$, либо в графе $\widetilde Y_{k-1}$, то Консерватор, очевидно, победит, выбрав соответствующую вершину либо в графе $\widetilde Y_{k-1}$, либо в графе $\widetilde X_{k-1}$.
Если же в $k$-ом раунде Новатор выберет вершину, не принадлежащую ни графу $\widetilde X_{k-1}$, ни графу $\widetilde Y_{k-1}$, то пусть для определенности он выбрал вершину в графе $X_{k-1}$. Тогда, в силу победы Консерватора в $(k-2)$-ом раунде, пара $\left(X_{k-1}|_{V(\widetilde
X_{k-1})\cup\{x_k^k\}},\widetilde X_{k-1}\right)$ является $(\alpha+\varepsilon)$--надежной. Следовательно, в силу свойства {\sf 2)}, в графе $Y_{k-1}$ существует такая вершина $y_k^k$, что граф $Y_{k-1}|_{V(\widetilde
Y_{k-1})\cup\{y_k^k\}}$ является точным $\left(X_{k-1}|_{V(\widetilde
X_{k-1})\cup\{x_k^k\}},\widetilde X_{k-1}\right)$--расширением графа $\widetilde Y_{k-1}$. Следовательно, Консерватор победит, выбрав вершину $y_k^k$.

Пусть в $(k-1)$-ом раунде Новатор выбирает вершину, принадлежащую либо
графу $X_{k-2}\setminus\widetilde X_{k-2}$ (положим $\widetilde X_{k-1}=\widetilde X_{k-2}$, $\widetilde Y_{k-1}=\widetilde Y_{k-2}$), либо графу $Y_{k-2}\setminus\widetilde Y_{k-2}$ (положим $\widetilde Y_{k-1}=\widetilde X_{k-2}$, $\widetilde X_{k-1}=\widetilde Y_{k-2}$). Для каждого $[s/t]$-элементного подмножества множества $\{x^{1}_{k-1}, x^{2}_{k-1},\ldots,x^{k-2}_{k-1}\}$ выберем в графе $X_{k-1}\setminus \widetilde{X}_{k-1}$ ровно одну вершину $z$ (если она существует), соединенную ребром со всеми вершинами из данного $[s/t]$-элементного подмножества и с вершиной $x_{k-1}^{k-1}$. Объединим все выбранные таким образом вершины вместе с вершиной $x_{k-1}^{k-1}$ в множество $V$. Так как $|V|\leq C^{[s/t]}_{k-2}+1=m$, то в силу определения графа $\widetilde X_{k-1}$ пара $\left(X_{k-1}|_{V(\widetilde X_{k-1})\cup V},\widetilde X_{k-1}\right)$ является $(\alpha+\varepsilon)$--надежной. Следовательно, в силу свойства {\sf 2)} в графе $Y_{k-1}$ существует такая вершина $y_{k-1}^{k-1}$ и граф $W$, являющийся точным $\left(X_{k-1}|_{V(\widetilde X_{k-1})\cup V},\widetilde X_{k-1}\right)$--расширением графа $\widetilde Y_{k-1}$, что при изоморфизме графов $X_{k-1}|_{V(\widetilde X_{k-1})\cup V}$ и $W$ вершины $x_{k-1}^1,x_{k-1}^2,\ldots,x_{k-1}^{k-1}$ переходят в вершины $y_{k-1}^1,y_{k-1}^2,\ldots,y_{k-1}^{k-1}$ соответственно.

Если в $k$-ом раунде Новатор выберет вершину либо в графе
$\widetilde X_{k-1}$, либо в графе $\widetilde Y_{k-1}$, то Консерватор, очевидно, победит, выбрав нужную вершину либо в графе $\widetilde Y_{k-1}$, либо в графе $\widetilde X_{k-1}$ соответственно. Если же в $k$-ом раунде Новатор выберет вершину, не принадлежащую ни графу $\widetilde X_{k-1}$, ни графу $\widetilde Y_{k-1}$, то пусть для определенности он выбрал вершину в графе $X_{k-1}$. Тогда в силу определения графа $\widetilde X_{k-1}$ вершина $x^{k}_{k}$ соединена ребром с не более, чем $[s/t]$ вершинами из графа $\widetilde X_{k-1}$, следовательно, она соединена ребром либо с не более, чем $[s/t]$ вершинами из множества $\{x^{1}_{k-1}, x^2_{k-1},\ldots,x^{k-1}_{k-1}\}$ либо с вершиной $x^{k-1}_{k-1}$ и с $[s/t]$ вершинами из множества $\{x^{1}_{k-1}, x^2_{k-1},\ldots,x^{k-2}_{k-1}\}$ (иначе
$$
 f_{\alpha-\varepsilon}(X_{k-1}|_{V(\widetilde X_{k-1})\cup\{x_k^k\}},\widetilde X_{k-1})\leq 1-(\alpha-\varepsilon)([s/t]+1)<1-\frac{([s/t]+1)tq}{sq+1}\leq1-\frac{(s/t+1/t)tq}{sq+1} <0,
$$
что противоречит определению графа $\widetilde X_{k-1}$). Во втором случае Консерватор сможет выбрать
нужную вершину $y_k^k$ и выиграть по построению. В первом же случае в силу свойства {\sf 2)} в графе $Y_{k-1}$ существует такая
вершина $y_k^k$, что граф $Y_{k-1}|_{\{y_{k}^1,y_{k}^2,\ldots,y_k^{k}\}}$ является точным $\left(X_{k-1}|_{\{x_{k}^1,x_{k}^2,\ldots,x_{k}^{k}\}},X_{k-1}|_{\{x_{k}^1,x_{k}^2,\ldots,x^{k-1}_{k}\}}\right)$--расширением графа $Y_{k-1}|_{\{y_{k}^1,y_{k}^2,\ldots,y_{k}^{k-1}\}}$. Следовательно, Консерватор победит, выбрав вершину $y_{k}^{k}$.

\section{Доказательство теоремы~\ref{refutation}}
\label{proof_refutation}

Пусть $G$ --- произвольный граф, $u,v\in V(G)$, $u\neq v$. Будем говорить, что $G$ является {\it $m$--цепью с концами $u$ и $v$}, если выполнено следующее условие. Для некоторого $d\in\mathbb{N}$ существуют такие подграфы $W,c_1,c_2,\ldots,c_n\subset G$, что $W$ --- цепь длины $d$ ({\it длиной} цепи называется количество ребер в ней) с концами $u,v$ и ребрами $e_1,e_2,\ldots,e_d$, для любого $i\in\{1,2,\ldots,d\}$ $c_i$ --- клика (полный граф) на $m$ вершинах, содержащая ребро $e_i$, и, кроме того, $G=c_1\cup c_2\cup\ldots\cup c_d$, любые две клики $c_i,c_j$, $i,j\in\{1,\ldots,d\}$, $i\neq j$, пересекаются по не более чем одной вершине. Наименьшее из всех чисел $d$, для которых выполнены эти условия будем называть {\it длиной} $m$--цепи $G$. Будем называть $m$--цепь длины $d$ {\it простой}, если в качестве $W$ можно взять простую цепь длины $d$. Заметим, что если $m$--цепь длины $d$ является простой, то любая клика $c_i$, $i\in\{1,\ldots,d\}$, из ее определения имеет общие вершины лишь с ``соседними'' кликами $c_{i-1}$ и $c_{i+1}$ (здесь и далее мы считаем, что $c_0=c_d$, $c_{d+1}=c_1$). Заметим, наконец, что у любой $m$--цепи $G$ с концами $u,v$ существует индуцированный подграф, который является простой $m$--цепью с концами $u,v$ (чтобы доказать существование такого подграфа достаточно в качестве графа $W$, определяющего этот подграф, выбрать цепь наименьшей длины с концами $u,v$ в $G$).

Будем говорить, что граф $G$ является {\it $m$--циклом}, если выполнено следующее условие. Для некоторого $d\in\mathbb{N}$ существуют такие подграфы $C,c_1,c_2,\ldots,c_d\subset G$, что $C$ --- простой цикл длины $d$ ({\it длиной} простого цикла называется количество ребер в нем) с ребрами $e_1,e_2,\ldots,e_d$, для любого $i\in\{1,2,\ldots,d\}$ $c_i$ --- клика на $m$ вершинах, содержащая ребро $e_i$, и, кроме того, $G=c_1\cup c_2\cup\ldots\cup c_d$, любая клика $c_i$, $i\in\{1,\ldots,d\}$, пересекается только с ``соседними'' кликами $c_{i-1},c_{i+1}$ и только по одной вершине (здесь и далее мы считаем, что $c_0=c_d$, $c_{d+1}=c_1$). Число $d$, для которого выполнены эти условия, будем называть {\it длиной} $m$--цикла $G$. Вершины графа $C$ в определении $m$--цикла мы будем называеть его {\it узловыми вершинами}. 
Заметим, что число вершин в $m$--цикле длины $d$ есть $d(m-1)$. Мы будем обозначать $M_d$ граф, являющийся $m$--циклом длины $d$.

Обозначим $k_1=k-10m+8$. Из условия теоремы имеем $k_1\geq 3$. Пусть $l_1,l_2\in[4,2^{k_1}]$ --- два натуральных числа. Рассмотрим граф $G_{l_1,l_2}$, являющийся объединением двух копий $M$ и $M'$ графов $M_{l_1},M_{l_2}$ соответственно, пересекающихся ровно по одной вершине, которая является узловой вершиной обоих графов. Плотность графа $G$ равна
$$
\dfrac{|E(G)|}{|V(G)|}=\dfrac{(l_1+l_2)\cdot\frac{m(m-1)}{2}}{(l_1+l_2)(m-1)-1}=
\frac{1}{\frac{2}{m}-\frac{2}{(l_1+l_2)m(m-1)}}\geq\frac{1}{\alpha},
$$
причем равенство достигается тогда и только тогда, когда $l_1=l_2=2^{k_1}$. Докажем, что граф $G$ является строго сбалансированным. Обозначим подграфы $C,c_1,\ldots$ из определения $m$--циклов $M$ и $M'$ за $C_1,c_1^1,\ldots,c_{l_1}^1$ и $C_2,c_1^2,\ldots,c_{l_2}^2$ соответственно. Предположим, что $G$ не строго сбалансирован. Тогда существует такой подграф $H\subset G$, что $\rho(H)\geq\rho(G)$. Можно считать, что $H$ --- связен, так как его плотность не превосходит максимальной из плотностей его компонент связности. Можно также считать, что для каждой из клик $c_i^j$, $i\in\{1,2,\ldots,l_j\}$, $j\in\{1,2\}$, хотя бы одна вершина которой входит в $H$, либо все ее ребра входят в $H$, либо все они в $H$ не входят, иначе добавив ``недостающие'' ребра, получим граф, плотность которого больше, чем $\rho(G)$. Действительно, предположим, что клика $c_i^j$ содержится в $H$ ``частично'', а именно в $H$ входит $v$ ее вершин и $\frac{v(v-1)}{2}$ ребер, где $v\in [2,m-1]$. Тогда, дополнив граф $H$ оставшейся частью клики $c_i^j$, мы получим граф с плотностью 
$$
 \frac{e(H)+\frac{m(m-1)}{2}-\frac{v(v-1)}{2}}{v(H)+(m-v)}>\rho(G),
$$
так как 
$$
 \frac{\frac{m(m-1)}{2}-\frac{v(v-1)}{2}}{m-v}=\frac{m+v-1}{2}\geq\frac{m+1}{2}>\rho(G)
$$ 
при $m\geq 2$. Можно также считать, что если $H$ содержит ровно одну вершину некоторой клики $c_i^j$, то он содержит не более одной вершины соседних с ним клик $c_{i-1}^j$ и $c_{i+1}^j$, иначе, удалив соседнюю клику (без одной вершины клики $c_i^j$), мы уменьшим число вершин на $m-1$, а число ребер на $\frac{m(m-1)}{2}$, следовательно, плотность графа увеличится. Таким образом, можем считать, что $H$ --- это один из $m$--циклов $M$ или $M'$. Но в таком случае $\rho(H)=\frac{m}{2}<\rho(G)$. Получили противоречие.

Поскольку $G_{2^{k_1},2^{k_1}}$ --- строго сбалансированный граф с плотностью $\frac{1}{\alpha}$, по теореме~\ref{small_subgraphs} имеем
\begin{equation}
\lim\limits_{n\rightarrow\infty}{\sf P}_{n,p}(N_G>0)=1-e^{-1/a(G)}\notin\{0,1\}.
\label{not_0_1}
\end{equation}

Будем говорить, что два свойства $A,B$ {\it асимптотически эквивалентны}, если пределы $\lim\limits_{n\rightarrow\infty}{\sf P}_{n,p}(A)$ и $\lim\limits_{n\rightarrow\infty}{\sf P}_{n,p}(B)$ существуют и равны. Далее мы приведем пример свойства $A$ из класса $\mathcal{L}_k$, асимптотически эквивалентного свойству $\{N_G>0\}$.\\

Пусть формула $\NI(u_1,\ldots,u_h)$ выражает свойство, заключающееся в том, что вершины $u_1,\ldots,u_h$ попарно различны: $\NI(u_1,\ldots,u_h)=\left(\bigwedge_{1\leq i<j\leq h}(\neg(u_i=u_j))\right)$. Определим формулу $\K(x_1,x_2,\ldots,x_m)$, выражающую свойство ``вершины $x_1,x_2,\ldots,x_m$ образуют клику размера $m$'':
\begin{align*}
 \K(x_1,x_2,\ldots,x_m)=\left(\NI(x_1,x_2,\ldots,x_m)\wedge\left(\bigwedge_{1\leq i<j\leq m}(x_i\sim x_j)\right)\right).
\end{align*}
Далее рассмотрим формулу $\MK(x_1,x_2,\ldots,x_m)$, выражающую свойство ``вершины $x_1,x_2,$ $\ldots,x_m$ образуют клику размера $m$, и любая клика графа либо совпадает с этой, либо пересекает ее по не более чем одной вершине'': $\MK(x_1,x_2,\ldots,x_m)=$
\begin{align*}
 \left(\K(x_1,x_2,\ldots,x_m)\wedge\left(\forall y_1\forall y_2\ldots \forall y_m\,\,
 \left(\K(y_1,y_2,\ldots,y_m)\Rightarrow \phi(x_1,\ldots,x_m,y_1,\ldots,y_m)\right)\right)\right),
\end{align*}
где формула $\phi(x_1,\ldots,x_m,y_1,\ldots,y_m)=$
\begin{align*}
 \left(\left(\bigvee_{\sigma\in\Sigma}\bigwedge_{i=1}^n
 (x_i=y_{\sigma(i)})\right)\vee\left(\bigvee_{i=1}^m\bigwedge_{1\leq \hat i\leq m,\,\hat i\neq i}\bigwedge_{j=1}^m(x_{\hat i}\neq y_{j})\right)\right)
\end{align*}
выражает свойство ``множества $\{x_1,\ldots,x_m\}$, $\{y_1,\ldots,y_m\}$ либо совпадают, либо пересекаются по не более чем одной вершине'' (здесь $\Sigma$ --- множество всех перестановок последовательности $1,\ldots,m$). Определим, наконец, для произвольного натурального числа $l$ формулу $\D_{l}(x_1,x_2)$, из истинности которой следует существование $m$--цепи длины $l$ с концами $x_1$ и $x_2$, никакая вершина которого не проходит через заданное множество вершин $u_1,\ldots,u_h$: при $l>1$
$$
 \D_{l}(x_1,x_2,u_1,\ldots,u_h) = \left(\exists y\:\left(\D_{\lfloor\frac{l}{2}\rfloor}(x_1,y)\wedge \D_{\lceil\frac{l}{2}\rceil}(y,x_2)\wedge\NI(x_1,x_2,y,u_1,\ldots,u_h)\right)\right),
$$ 
$$
 \D_{1}(x_1,x_2,u_1,\ldots,u_h) = \left(\exists x_3\exists x_4\ldots\exists x_{m}\:\left(\MK(x_1,x_2,\ldots,x_m)
 \wedge\NI(x_1,\ldots,x_m,u_1,\ldots,u_h)\right)\right).
$$
Заметим, что если в графе существует $m$--цепь $G$ длины $l$ с концами $x_1$ и $x_2$, никакая вершина которого не проходит через заданное множество вершин $u_1,\ldots,u_h$, но не является истинной формула $\D_{l}(x_1,x_2,u_1,\ldots,u_h)$, то существует клика на $m$ вершинах, пересекающая граф $G$ по хотя бы одному ребру, но не содержащаяся в нем целиком.
Из определения формулы $\D_l(x_1,x_2,u_1,\ldots,u_h)$ 
легко видеть, что ее кванторная глубина 
равна $\lceil\log_2 l \rceil+2m-2$.


Определим, наконец, искомое свойство~$A$. Положим $l=2^{k_1-1}-2$.
$$
A=\Bigg(\exists x\exists y\exists y'
\exists v_1^1\ldots\exists v_1^{m-1}\ldots
\exists v_4^1\ldots\exists v_4^{m-1}
\exists u_1^1\ldots\exists u_1^{m-1}\ldots
\exists u_4^1\ldots\exists u_4^{m-1}
$$
$$
 \Bigg[\left(\NI(x,y,y',
 v_i^j,u_i^j,\,\,i\in\{1,2,3,4\},j\in\{1,\ldots,m-1\})\right)\wedge
$$
$$
 \left(\bigwedge_{i=1}^2\left(\MK(y,u_i^1,\ldots,u_i^{m-1})\wedge
 \MK(y',u_{i+2}^1,\ldots,u_{i+2}^{m-1})\right)\right)\wedge
 \left(\bigwedge_{i=1}^4 \left(\MK(x,v_i^1,\ldots,v_i^{m-1})\wedge\right.\right.
$$
$$
 \D_l\left(v_i^{m-1},u_i^{m-1},
 \left.\left.x,y,y',v_{\hat i}^j,u_{\hat i}^j,\,\,(\hat i,j)\in J\right)\right)\right)
 \Bigg]
 \Bigg),
$$
где $J=\{1,2,3,4\}\setminus\{i\}\times\{1,\ldots,m-1\}\cup\{i\}\times\{1,\ldots,m-2\}$.
Из определения свойства $A$ легко видеть, что его кванторная глубина равна $3+8m-8+\lceil\log_2 (2^{k_1-1}-2)\rceil+2m-2=10m-7+k_1-1=k_1+10m-8=k$. Следовательно, $A\in\mathcal{L}_{k}$.

Рассмотрим множество $\tilde\Omega_n\subset\Omega_n$ всех графов, в которых не существует подграфов $H$ с $v(H)\leq 2^{m(2k_1+1)}$ и $\rho(H)>1/\alpha$. В силу теоремы~\ref{erdos} имеем ${\sf P}_{n,p}(\tilde\Omega_n)=1$. Пусть $\mathcal{G}\in\tilde\Omega_n$ и $\mathcal{G}$ содержит подграф $X$, изоморфный графe $G_{2^{k_1},2^{k_1}}$. Пусть, кроме того, $\mathcal{K}$ --- множество всех таких пар $(K,T)$, что $K$ --- полный граф на $m$ вершинах, а $v(T)\geq 2$. Для того, чтобы доказать, что граф $\mathcal{G}$ обладает свойством $A$ достаточно доказать, что в $\mathcal{G}$ нет $(K,T)$--расширений подграфов графа $X$. Предположим, что хотя бы одно такое расширение $Y$ имеется. Тогда
$\rho(Y\cup X)>1/\alpha$, так как $\rho(X)=1/\alpha$, 
$$
 f_{\alpha}(Y,X)\leq v(Y,X)-\alpha\frac{v(Y,X)(2m-v(Y,X))}{2}\leq v(Y,X)\left(1-\alpha\frac{m+2}{2}\right)<0.
$$
Кроме того, $v(Y)\leq 2^{2k_1(m-1)-1+m-2}<2^{m(2k_1+1)}$
Тем самым, получили противоречие с принадлежностью графа $\mathcal{G}$ множеству $\tilde\Omega_n$, и, следовательно, граф $\mathcal{G}$ обладает свойством $A$.

Предположим теперь, что $\mathcal{G}\in\tilde\Omega_n$ и $\mathcal{G}$ обладает свойством $A$. Тогда в графе $\mathcal{G}$ найдутся такие попарно различные вершины $x,y,y',v_i^j,u_i^j,\,\,i\in\{1,2,3,4\},j\in\{1,\ldots,m-1\}$, что существуют $m$--цепи $C_i$, $i\in\{1,2,3,4\}$, длины $l$ с концами $v_i^{m-1},u_i^{m-1}$ (обозначим графы $W,c_1,\ldots,c_l$ из определения этих $m$--цепей $W_i,c_1^i,\ldots,c_l^i$), наборы вершин $x,v_i^1,\ldots,v_i^{m-1}$, $i\in\{1,2,3,4\}$, а также $y,u_i^1,\ldots,u_i^{m-1}$, $i\in\{1,2\}$, и $y',u_i^1,\ldots,u_i^{m-1}$, $i\in\{3,4\}$, образуют клики, причем графы $C_i$, $i\in\{1,2,3,4\}$, пересекают эти клики только по концевым вершинам цепей $W_i$. Как было замечено выше, в каждом из графов $C_i$, $i\in\{1,2,3,4\}$, содержится подграф $\tilde C_i$, являющийся простой $m$--цепью с концами $v_i^{m-1},u_i^{m-1}$. Если $\mathcal{G}$ не содержит подграфа, изоморфного графу $X$, то либо хотя бы одна из цепей $\tilde C_i$, $i\in\{1,2,3,4\}$, имеет длину меньше, чем $2^{k_1-1}-2$, либо хотя бы два графа из $\tilde C_i$, $i\in\{1,2,3,4\}$, пересекаются по вершинам. Очевидно, в таком случае в графе $\mathcal{G}$ существует подграф, изоморфный графу $G_{l_1,l_2}$, где $l_1,l_2\in[4,2^{k_1}]$, причем хотя бы одно из $l_1,l_2$ строго меньше чем $2^{k_1}$. Но $\rho(G_{l_1,l_2})>1/\alpha$, $v(G_{l_1,l_2})<2^{2k_1m}$. Получили противоречие. Следовательно, граф $\mathcal{G}$ содержит подграф, изоморфный $G_{2^{k_1},2^{k_1}}$.

Таким образом, в силу~(\ref{not_0_1})
$$
 {\sf P}_{n,p}(A)\sim {\sf P}_{n,p}(\tilde\Omega_n\cap A)={\sf P}_{n,p}(\tilde\Omega_n\cap \{N_G>0\})\sim
 {\sf P}_{n,p}(N_G>0)\sim 1-e^{-1/a(G)}.
$$
Следовательно, поскольку кванторная глубина свойства $A$ не превосходит $k$, то для $\alpha = \frac{2}{m} - \frac{1}{2^{k_1}m(m-1)}$ случайный граф $G(n,n^{-\alpha})$ не подчиняется $k$--закону нуля или единицы.

\end{document}